\documentclass{elsart}
\usepackage{ifpdf}
\usepackage{graphicx,natbib,amssymb,lineno}
\usepackage {amsmath, latexsym, amscd, amssymb}
\usepackage {graphicx}
\usepackage[latin1]{inputenc}
\usepackage [french]{babel}
\ifpdf
\usepackage[%
  pdftitle={Repr\'esentations lin\'eaires des graphes finis},%
  pdfauthor={Lucas Vienne},%
  pdfsubject={The preprint document class elsart},%
  pdfkeywords={automorphisme,  groupe, graphe, transitif},%
  pdfstartview=FitH,%
  bookmarks=true,%
  bookmarksopen=true,%
  breaklinks=true,%
  colorlinks=true,%
  linkcolor=blue,anchorcolor=blue,%
  citecolor=blue,filecolor=blue,%
  menucolor=blue,pagecolor=blue,%
  urlcolor=blue]{hyperref}
\else
\usepackage[%
  breaklinks=true,%
  colorlinks=true,%
  linkcolor=blue,anchorcolor=blue,%
  citecolor=blue,filecolor=blue,%
  menucolor=blue,pagecolor=blue,%
  urlcolor=blue]{hyperref}
\fi

\makeatletter
\def\elsartstyle{%
    \def\normalsize{\@setfontsize\normalsize\@xiipt{14.5}}
    \def\small{\@setfontsize\small\@xipt{13.6}}
    \let\footnotesize=\small
    \def\large{\@setfontsize\large\@xivpt{18}}
    \def\Large{\@setfontsize\Large\@xviipt{22}}
    \skip\@mpfootins = 18\p@ \@plus 2\p@
    \normalsize
}
\@ifundefined{square}{}{\let\Box\square}
\makeatother

\newcommand{\R}{\mathbb{R}}

\newcommand{\EE}{\mathcal{E}}

\newcommand{\GG}{\mathcal{G}}

\newcommand{\SSS}{\mathcal{S}}

\newcommand{\XX}{\mathcal{X}}

\newcommand{\ssi}{\Longleftrightarrow}
\newcommand{\x}{\times}

\newcommand{\al}{\langle}
\newcommand{\ar}{\rangle}

\newcommand{\Aut}{\mathrm{Aut}}

\newcommand{\stab}{\mathrm{stab}}
\newcommand{\Gram}{\mathrm{Gram}}

\newcommand{\pb}{\pagebreak}
\newcommand{\npb}{\nopagebreak}

\newcommand{\mm}{m\^eme }
\newcommand{\elt}{\'el\'ement }
\newcommand{\elts }{\'el\'ements }
\newcommand{\si }{si et seulement si }
\newcommand{\cad }{c'est-\`a-dire }

\newcommand{\Dem }{{\it D\'emonstration }}

\newcommand{\ttt}{th\'eor\`eme }

\newcommand{\oo}{\circ }
\newcommand{\cc }{\c{c}}

\newcommand{\ga }{\alpha }
\newcommand{\gb }{\beta }

\newcommand{\gd }{\delta }

\newcommand{\gve }{\varepsilon }
\newcommand{\gvf }{\varphi }

\newcommand{\go }{\omega }
\newcommand{\gs}{\sigma}

\newcommand{\gC }{\Gamma }

\newtheorem{deff}{D\'efinition}

\newtheorem{thmf}{Th\'eor\`eme }

\newtheorem{propf}{Proposition}

\begin{document}

\begin{frontmatter}
\title{Repr\'esentations lin\'eaires des graphes finis}

\author{Lucas Vienne}
\address{Departement de math\'ematiques. Universit\'e d'Angers. France}

\ead{lucas.vienne@univ-angers.fr}

\begin{abstract}
Let $X$ be a non-empty finite set and $\ga$ a symmetric bilinear form on a real finite dimensional vector space $E$. We say that a set $\GG=\{U_i \  |\  i\in
X\}$  of linear lines in  $E$ is an isometric sheaf, if  there exist generators $u_i$  of the lines $U_i$, and real constants $\go$  and $c$ such that \\[0,1cm]
\vspace{0,2cm}\centerline{$\forall i,j \in X,$ \  $\ga(u_i,u_i)=\go $, and  if \ $i\neq j$  \ then  \  $\ga(u_i,u_j)=\gve_{i,j}.c,$  \ with \  $\gve_{i,j}\in \{-1,+1\}$  }

  Let $\gC$ be the graph whose set of vertices is $X$,  two of them, say  $i$ and $j$, being linked when  $\gve_{i,j} =- 1$.
   In this article we explore the relationship between $\GG$ and $\gC$
   ; we describe all sheaves associated with a given graph  $\gC$ and
   construct the group of isometries stabilizing one of those as an
   extension  group of $\Aut(\gC)$. We finally illustrate our
   construction with some examples.
 \end{abstract}

\begin{keyword}
groupe, graphe, automorphisme, transitif
\end{keyword}
\end{frontmatter}
\section{Introduction}

Un entier $n\geq 1$ \'etant choisi, on pose  $X=\{1,\ldots,n\}$ et on
d\'efinit un graphe $(\gC,X)$  par l'ensemble $\gC$ de ses ar\`etes qui sont
des paires $\{x,y\}$  de points distincts pris dans l'ensemble
 $X$ de ses sommets. On notera parfois plus simplement $\gC$ le graphe
 $(\gC,X)$. Sa matrice $\EE$ est donn\'ee par son $(i,j)$-\`eme coefficient  $\gve_{i,j}$ qui
vaut\ $-1$ \  si $i$ et $j$ sont li\'es (not\'e  $i\sim j$) et $1$
dans le cas contraire.\\
Donnons nous  un espace quadratique $(E,\ga)$ \cad un espace
vectoriel  r\'eel $E$ de dimension finie muni d'une forme bilin\'eaire
sym\'etrique $\ga $.
\begin{deff}[Repr\'esentation d'un graphe]$\ $\\
 Soient $\go $ et $c$ deux constantes r\'eelles et  $(E,\ga)$ un
 espace quadratique. Une  application $u$ de $X$ dans
 $(E,\ga)$ not\'ee $u: i\to u_i$ est une {\it repr\'esentation} du
 graphe  $(\gC,X)$ de param\`etres  $(\go, c)$  si \\[0,1cm]
$(1)$
\vspace{0,2cm}\centerline{$\forall i,j \in X,$ \  $\ga(u_i,u_i)=\go $
   \ et \  si \ $i\neq j$,   \  $\ga(u_i,u_j)=\gve_{i,j}.c.$ }
La matrice $S(u)$ de coefficient g\'en\'eral $\ga(u_i,u_j)$ $(i,j\in X)$ est
dite {\it matrice de la repr\'esentation} $u$. On appellera respectivement {\it rang} et  {\it degr\'e} de $u$  le rang  de $S(u)$ et  la dimension de $E$.
\end{deff}

On associe \`a chaque repr\'esentation  $u :X\to E$ d'un graphe $(\gC,X)$ la gerbe des droites $U_i$ engendr\'ees par les $u_i$ ($i\in X$) :

\begin{deff}[Gerbe isom\'etrique]$\ $\\
Un ensemble $\GG=\{U_i\ |\ i\in X\}$ de droites vectorielles index\'e
par $X$ est une gerbe isom\'etrique, ou une $\ga$-gerbe dans l'espace
quadratique $(E,\ga)$,  s'il existe une repr\'esentation $u$ d'un graphe
$(\gC,X)$ dans l'espace $(E,\ga)$ telle que pour tout indice $i\in X$,
$U_i=\al u_i\ar$. On note alors $\GG(u)$ la gerbe $\GG$ que l'on dit
associ\'ee \`a la repr\'esentation $u$ de $\gC$.
\end{deff}

Dans cet article  on d\'ecrit toutes les gerbes associ\'ees \`a un
graphe $\gC$ donn\'e. On montre aussi que le groupe $H(\gC)$ des automorphismes d'un graphe
$(\gC,X)$ poss\`ede une image naturelle $H(u)$ dans le groupe
d'automorphismes $\Aut(\GG(u))$ de chaque gerbe $\GG(u)=\{\al u_1
\ar,\ldots,\al u_n\ar\}$ qui lui est associ\'ee. Le \ttt 2 prouve qu'en g\'en\'eral le groupe $\Aut(\GG(u))$ est 
beaucoup plus gros que $H(u)$ et que sa structure ne d\'epend
pas de la repr\'esentation $u$, sauf lorsque les param\`etres $\go$ et
$c$ de $u$ ont le \mm module ($|\go|=|c|$).
Pour finir, quelques exemples illustrent l'extension du groupe $H(u)$
au groupe $\Aut(\GG(u))$ ; si $\gC$ est un carr\'e, le groupe $\Aut(\GG(u))$ s'identifie
au groupe d'isom\'etries d'un cube tandis que si $\gC$ est un
pentagone, le groupe $\Aut(\GG(u))$ s'identifie
au groupe d'isom\'etries d'un dod\'eca\`edre.\\
Un deuxi\`eme article, en pr\'eparation, classe les graphes conduisant \`a un groupe $\Aut(\GG(u))$ op\'erant doublement transitivement sur la gerbe $\GG(u)$.

\subsection{D\'etermination des repr\'esentations d'un graphe $(\gC,X)$} 
{\it Terminologie.} 
Une forme bilin\'eaire sym\'etrique $\ga$ \'etant donn\'ee sur
un espace vectoriel r\'eel $E$, on appelle {\it matrice de Gram } d'un $n$-uplet de vecteurs
$(u_1,\ldots, u_n)$ la matrice $\Gram_\ga(u_1,\ldots, u_n)$  dont le coefficient de position
$(i,j)$ est $\ga(u_i,u_j)$. Le r\'esultat suivant est obtenu
par des arguments \'el\'ementaires :

\begin{propf}$\ $\\ 
Soit $S=(S_{i,j})_{1\leq i,j\leq n}$ une matrice sym\'etrique de rang
$r$ dans $M_n(\R)$ et $E=\R^r$.\\
$1$. Il existe une forme bilin\'eaire sym\'etrique $\ga$ sur $E$ et
des vecteurs $u_1,\ldots, u_n$ dans $E$ tels que \quad $S=\Gram_\ga(u_1,\ldots, u_n)$\\
  $2$. Si $S$ est aussi la matrice de Gram d'un autre syst\`eme de
vecteurs $v_1,\ldots, v_n$ de $E$ pour une forme bilin\'eaire $\beta $, il  existe une isom\'etrie $f : (E,\ga) \to
(E,\gb)$ telle que  pour tout indice $i$ dans $ X=\{1,\ldots,n\}$, \quad  $f(u_i)=v_i$
 \end{propf}

\subsubsection{Notion d'isomorphisme. Op\'erations sur les repr\'esentations} 

\begin{deff}$\ $\\ 
Deux repr\'esentations $u : X\to (U,\ga)$ et $v :X\to (V,\gb)$ d'un graphe
$(\gC,X)$ dans des espaces quadratiques $(U,\ga)$ et $(V,\gb)$ sont dites isomorphes s'il existe une isom\'etrie $f $ de $(U,\ga)$ sur $(V,\gb)$ telle que \ $f\oo u=v$.
\end{deff}

{\it Repr\'esentation triviale, repr\'esentation nulle} \\
Soit $u :X \to (E,\ga)$ une repr\'esentation d'un graphe $(\gC,X)$ de
param\`etres $(\go,c)$. Dire que les vecteurs $u_i$ ($i\in X$) sont
deux \`a deux orthogonaux \'equivaut \`a dire que le param\`etre $c$
est nul. Dans ce cas nous dirons que la repr\'esentation $u$ est {\it triviale}. 
 Si maintenant la forme bilin\'eaire $\ga$ est nulle, les deux
 param\`etres $\go$ et $c$ de $u$ sont nuls, et nous dirons que la
 repr\'esentation $u$ est {\it nulle}.
 
{\it Somme} \\
Soient   $u : X \to (U,\ga)$ et  $v :X\to (V,\gb)$ deux repr\'esentations
d'un graphe  $(\gC,X)$ associ\'ees aux param\`etres $(\go _u, c_u)$ et
$(\go _v,c_v)$, et $U\oplus V$ la somme directe orthogonale de
$U$ et $V$. L'application  $w : X \to U\oplus V$  donn\'ee par
$w_i=u_i+v_i$ pour tout $i\in X$ est une repr\'esentation de $(\gC,X)$
dans l'espace $U\oplus V$ de matrice $S(w)=S(u)+S(v)$ et de param\`etres $\go _w=\go _u+\go _v$
et $c_w=c_u+c_v$. On dit que $w$ est la {\it somme} de $u$ et $v$, et on la  note $w=u+v$.

{\it Remarque.} L'addition des repr\'esentations d'un graphe $\gC$ est
associative, \`a isomorphisme pr\`es : si $u,v$ et $w$ sont trois
repr\'esentations de $\gC$ alors \\[0.1cm]
\centerline{$(u+v)+w\simeq u+(v+w)$.}

{\it Repr\'esentation r\'eduite} \\
Soit $u :X\to (E,\ga)$ une repr\'esentation de rang $r$ d'un graphe
$(\gC,X)$, $V$ un suppl\'ementaire dans $E$ de l'espace $U=\al
u_1,\ldots,u_n\ar$, $N=U\cap U^\perp$ et $R$ un suppl\'ementaire de $N$
dans $U$ qui est donc de dimension $r$.
 Notant $p$ et $q$ les projections sur $R$
et  $N \oplus V$ associ\'ees \`a la d\'ecomposition $E=R\oplus( N \oplus V)$, on v\'erifie
simplement que  \\
$*$ \  $p\oo u : X \to R$ \ est une repr\'esentation de $\gC$ de matrice
$S(p\oo u)=S(u)$ et de rang $r=\dim R$.\\
$*$ \  $q\oo u : X \to  N \oplus V$ \  est une repr\'esentation nulle.\\
$*$ \  $u=p\oo u + q\oo u$ est 
donc $u$ est la somme d'une repr\'esentation de degr\'e $r$ et d'une
repr\'esentation nulle. De plus, pour toute autre d\'ecomposition de
$U$ en somme $R'\oplus N$ les repr\'esentations de $\gC$ sur $R$ et
$R'$ sont visiblement isomorphes. R\'esumons  :

\begin{propf}[Repr\'esentation r\'eduite]$\ $\\  
Une repr\'esentation $u : X \to (E,\ga)$ est dite {\bf r\'eduite} si son
degr\'e $d=\dim E$ est \'egal au rang $r$  de sa matrice $S(u)$.
Toute repr\'esentation $u$ d'un graphe $(\gC,X)$ est la somme d'une  repr\'esentation nulle  et d'une repr\'esentation  r\'eduite, uniquement d\'etermin\'ee \`a isomorphisme pr\`es.
\end{propf}
 
 \subsubsection{Classification} 

Soit $(\gC,X)$ un graphe et $u :X\to (E,\ga)$ une repr\'esentation
r\'eduite de $\gC$. Sa matrice $S(u)$ ne d\'ependant que des
param\`etres $\go $ et $c$  (d'apr\`es $(1)$), on la note aussi
$S(u)=S(\go,c)$. Donnons nous inversement un couple de r\'eels
$(\go,c)$ et la matrice $S(\go,c)$ qui lui est associ\'ee. La
proposition 1 nous montre qu'\`a isomorphisme pr\`es, il existe une
unique repr\'esentation r\'eduite  $u :X\to (E,\ga)$ dont la matrice
est $S(u)=S(\go,c)$. Ainsi,

\begin{thmf}$\ $\\ 
 Soit  $(\gC,X)$ un graphe. \\
 $1.$  L'application qui associe
 \`a toute repr\'esentation $u : X\to E$ de $\gC$ le couple $(\go ,c)$
 de ses param\`etres dans $\R^2$ est surjective. Deux repr\'esentations r\'eduites de $\gC$ sont isomorphes \si elles ont
 les m\^emes param\`etres.\\
 $2.$ Le degr\'e d'une repr\'esentation r\'eduite de matrice $S(1,c)$ est le rang de cette matrice, \cad  $n-\mu(c)$ o\`u $n=|X|$ et $\mu(c)$ est la multiplicit\'e de $c$ comme racine du polyn\^ome $\chi(x)=\det(S(1,x)$.
  \end{thmf}

\subsection{Automorphismes d'un graphe et des gerbes associ\'ees}

{\it Notations usuelles. } \`A toute permutation $\gs \in \SSS_X$
associons sa matrice $P_\gs$ qui est de type $n\x n$ et dont le $(j,i)$-\`eme terme vaut
$P_{\gs,j,i}=\gd_{j,\gs(i)}$ (o\`u $\gd$ est la fonction de Kronecker
usuelle). Pour toute matrice $M$ de type $n\x n$ posons  \ ${^\gs
  M}=P_\gs.M.P_\gs^{-1}$, et notons enfin $\stab(M)$ le stabilisateur de $M$ dans $\SSS_X$, \cad
l'ensemble des permutations $\gs\in \SSS_X$ telles que \  ${^\gs M}=M$.\\[-0,8cm]

 \subsubsection{Automorphismes d'un graphe $(\gC,X)$}
 Soit  $\EE$ la matrice d'un graphe $\gC$ et $H(\gC)$ le
groupe de ses {\it automorphismes}, \cad le sous-groupe du
groupe $\SSS_X$ des  permutations de $X=\{1,\ldots,n\}$ qui pr\'eservent l'ensemble $\gC$ de ses ar\`etes.
On v\'erifie que pour $\gs$ dans $\SSS_X$, la matrice $^\gs \EE$ est la matrice  du graphe $\gs(\gC)$, donc le stabilisateur dans $\SSS_X$ du graphe  $\gC$ est aussi le stabilisateur de sa matrice $\EE$ : \\[0,1cm]
\vspace{-0,5cm} \centerline{$H(\gC)=\stab(\EE)$, \  ou encore : \ $\forall 
    \gs \in \SSS_X,\    {^\gs\EE}=\EE \ssi  \gs(\gC)=\gC$.}
    
  \subsubsection{Automorphismes d'une gerbe} 
  
 Soit $u :X\to (E,\ga)$ une repr\'esentation r\'eduite non triviale
 d'un graphe $(\gC,X)$ o\`u, pour \'eviter des discussions sans grand int\'er\^et, on suppose que $|X|\geq 3$.\\
 On appelle {\it automorphisme} de la gerbe  $\GG(u)=\{\al
 u_1\ar\ldots,\al u_n \ar\}$ associ\'ee \`a la  repr\'esentation $u$
 toute isom\'etrie $\gvf$ de $(E,\ga)$ qui induit une permutation de l'ensemble
 $\GG(u)$ par $\gvf(\al u_i \ar)=\al \gvf(u_i) \ar$ (pour tout $i\in X$). On note $\Aut(\GG(u))$ le groupe de ces automorphismes.\\
 
{\it Un exemple.}  
Choisissons une permutation $\gs$ dans le groupe $H=H(\gC)$ des
automorphismes du graphe $(\gC,X)$. Les matrices $\EE$
et $^\gs\EE$ \'etant \'egales, les
syst\`emes de vecteurs $(u_1,\ldots,u_n)$ et
$(u_{\gs(1)},\ldots,u_{\gs (n)})$ ont \mm matrice de Gram et \mm rang, $r=\dim E$. Il existe donc, d'apr\`es la proposition 1, une
isom\'etrie $f_\gs$ de $(E,\ga)$ v\'erifiant, pour tout $i\in X$,
$f_\gs(u_i)=u_{\gs(i)}$. Par suite $f_\gs$  est un automorphisme
de la gerbe $\GG(u)$, et on v\'erifie simplement que l'application $f :H\to \Aut(\GG(u)) $ qui associe
$f_\gs$ \`a chaque $\gs\in H$ est un homomorphisme de groupes.

{\it Une g\'en\'eralisation.} 
\'Etudions les automorphismes $\gvf$ de la gerbe  $\GG(u)$ tels que, pour une permutation convenablement choisie  $\gs \in \SSS_X$, on ait \\[0.1cm]
$(3)$ \vspace{0,1cm}\centerline{ $\forall i \in X$, \quad
  $\al\gvf(u_i)\ar=\al u_{\gs(i)}\ar$.}
Il existe, dans ce cas, un syst\`eme de scalaires  $\nu=(\nu_1,\ldots,\nu_n)$ tels que \\[0.1cm]
$(4)$ \vspace{0,1cm}\centerline{ $\forall i \in X$, \quad $\gvf(u_i
)=\nu_i.u_{\gs(i)}$.}
Montrons qu'alors les $\nu_i$, pour $i\in X$,  sont tous dans l'ensemble $\{-1,1\}$. D'apr\`es la \mbox{proposition 1,}  $\gvf$ est une
isom\'etrie, donc un automorphisme de la gerbe $\GG(u)$, \si 
 les matrices de Gram des syst\`emes de vecteurs $(u_1,\ldots,u_n)$ et
$(\nu_1.u_{\gs(1)},\ldots,\nu_n.u_{\gs(n)})$ sont \'egales, autrement
dit si pour $i\neq j$, \\[0,2cm]
\vspace{0,2cm}\centerline{ $\ga(u_i,u_j)=\ga(f(u_i),f(u_j))=\ga(\nu_i.u_{\gs(i)},\nu_j.u_{\gs(j)})=\nu_i.\nu_j.\gve_{\gs(i),\gs(j)}.c=\gve_{i,j}.c$,}
soit encore puisque, la repr\'esentation $u$ \'etant non triviale,
$c$ est non nul : \\[0.2cm]
$(5)$ \vspace{0.2cm}\centerline{  $\forall i,j \in X$, \quad
  $\nu_i.\nu_j.\gve_{\gs(i),\gs(j)}=\gve_{i,j}$.}
Mais comme $|X|\geq 3$, pour trois indices $i,j,k$ distincts dans $X$, il vient
$1=|\nu_i.\nu_j|=|\nu_i.\nu_k|=|\nu_j.\nu_k|$,  d'o\`u l'on tire \ 
$1=|\nu_i|=|\nu_j|=|\nu_k|$, donc le $n$-uplet
$\nu=(\nu_1,\ldots,\nu_n)$ est \`a valeurs dans $\{-1,1\}$. 
La condition $(5)$ nous sugg\`ere la 

\begin{deff}$\ $\\ 
Deux matrices $M$ et $N$ de \mm type $n\x n$ sont dites associ\'ees
s'il existe une suite $\nu_1,\ldots,\nu_n$ de
 coefficients, tous pris dans $\{-1,1\}$ tels que \\[0.2cm]
\centerline{$\forall i, j ,\   1\leq i,j\leq n, \qquad M_{i,j}=\nu_i.\nu_jN_{i,j}$}
\end{deff}
On peut alors r\'esumer la discussion qui pr\'ec\`ede par la 
\begin{propf}$\ $\\ 
Un automorphisme $\gvf$ de la gerbe $\GG(u)$ est associ\'e \` a une
permutation $\gs $ dans $\SSS_X$ par la relation $(3)$ \si les matrices $\EE$ et
$^\gs\EE $ sont associ\'ees. Il existe, dans ce cas, un unique $n$-uplet
$\nu=(\nu_1,\ldots,\nu_n)$ de coefficients tous pris dans $\{-1,1\}$
tels que \\[-0.1cm]
\vspace{0,1cm}\centerline{ $\forall i \in X$,\quad  $\gvf(u_i
)=\nu_i.u_{\gs(i)}$.}
L'isom\'etrie $\gvf$ est dite {\it associ\'ee} au couple
$(\gs,\nu)$ et  not\'ee $\gvf=f_{\gs,\nu}$.
\end{propf}
On munit maintenant l'ensemble  $G$  des couples $(\gs,\nu)\in \SSS_X\x
\{-1,1\}^n$ satisfaisant \`a la condition $(5)$ d'une structure de
groupe pour laquelle l'application $f :G\to \Aut(\GG(u))$ donn\'ee par
$f(\gs,\nu)=f_{\gs,\nu}$ est un morphisme de groupes. Choisissant deux \elts $(\gs,\nu)$ et $(\gs',\nu')$
  dans $G$, les isom\'etries $f_{\gs,\nu}$ et $f_{\gs',\nu'}$ de
  $(E,\ga)$ qui leurs sont associ\'ees v\'erifient \\[0.1cm]
\vspace{0,1cm}\centerline{$f_{\gs,\nu}\oo
  f_{\gs',\nu'}=f_{\gs'',\nu''}$ \ avec \ $\gs''=\gs\oo\gs'$ \ et pour
  tout $ i\in X,\ \nu''_i=\nu_{\gs'(i)}.\nu'_i$.}
Ceci nous conduit \`a noter  $\nu^{\gs'} $ le $n$-uplet donn\'e par\\[0.1cm]
\vspace{0,1cm}\centerline{$\forall i\in X$, \quad
  $\nu^{\gs'}_i=\nu_{\gs'(i)} $, }
et \`a d\'efinir sur $G$ la loi $*$  en posant, pour $(\gs,\nu)$ et $(\gs',\nu')$  dans  $G$\\[0.2cm]
$(6)$\vspace{0,2cm}\centerline{$ (\gs,\nu)*(\gs',\nu')=(\gs\oo\gs',
  \nu^{\gs'}.\nu') $.}
On v\'erifie alors simplement la 
\begin{propf}$\ $\\ 
L'ensemble  $G$  des couples $(\gs,\nu)\in \SSS_X\x
\{-1,1\}^n$ satisfaisant \`a la condition $(5)$ est mumi d'une
structure de groupe par la relation $(6)$, pour laquelle l'application $f:(\gs, \nu)\to f_{\gs,\nu}$ devient
un morphisme de $G$ dans le groupe $\Aut(\GG(u))$.\\
On note $G(u)$ l'image de $G$ dans $\Aut(\GG(u))$ par ce morphisme
 \end{propf} 
Nous pouvons maintenant compl\'eter la description du groupe
$\Aut(\GG(u))$ :
\begin{thmf}$\ $\\ 
Soit $u $ une repr\'esentation  r\'eduite et non triviale d'un  graphe
$(\gC,X) $ dans un espace quadratique $(E,\ga) $. On suppose de plus
$|X|\geq 3$. \\
$1.$ Lorsque les droites  $\al u_i\ar$ de la gerbe $\GG(u)=\{\al u_1\ar\ldots,\al u_n \ar\}$ sont deux \`a deux distinctes,  le morphisme $f :G\to \Aut(\GG(u)) $ est un isomorphisme  et on a \\[0,1cm]
 \vspace{0,1cm}\centerline{$G\simeq G(u) = \Aut(\GG(u))$.}
 $2.$ Lorsque les modules des param\`etres $\go$ et $c$ de $u$ sont
distincts $(|\go|\neq |c|)$, les droites de la gerbe $\GG(u)$ sont deux \`a deux distinctes.
\end{thmf}
\Dem \\
$1.$  Par d\'efinition, tout automorphisme $\gvf$ de la gerbe
$\GG(u)=\{\al u_1\ar\ldots,\al u_n \ar\}$ induit une permutation
de l'ensemble $\GG(u)$. Lorsque les $n$ droites $\al u_1\ar\ldots,\al
u_n \ar$ sont deux \`a deux distinctes il existe donc une permutation
$\gs$ de $X$ telle que pour tout indice $i$, $\gvf(\al u_i \ar)=\al
u_{\gs(i)}\ar$, ce qui n'est rien d'autre que la relation $(3)$
ci-dessus. La proposition 3 nous montre qu'alors il existe un \elt
$(\gs,\nu)$ dans $G$ tel que $\gvf=f_{\gs,\nu}$, et ceci prouve la
surjectivit\'e du morphisme $f$.\\
Un \elt $(\gs,\nu)$ du groupe $G$ est dans le noyau du morphisme
$f$ \si pour tout indice $i\in X$,\  $f_{\gs,\nu}(u_i)=\nu_i.u_{\gs(i)}= u_i$. Mais les $n$ droites $\al u_i\ar$  \'etant distinctes ces
\'egalit\'es impliquent  $\gs(i)=i $ et $\nu_i=1$ pour tout indice $i$
dans $X$. D'o\`u l'injectivit\'e.\\
$2.$ On montre maintenant, par contraposition, que si $(|\go|\neq |c|)$, les
droites $\al u_1\ar\ldots,\al u_n \ar$ sont deux \`a deux distinctes.
Si, pour deux indices $i$ et $j$ distincts, on a  $\al u_i\ar= \al u_j\ar$,
il existe un nombre \ $\gve$ non nul  tel que \ $u_j=\gve u_i$ \  et donc\\[0.1cm]
\vspace{0,1cm}\centerline{$\ga(u_i,u_j)=\gve.\ga(u_i,u_i)=\gve.\go=\gve_{i,j}.c$, 
  \ o\`u \ $\gve_{i,j}=\pm 1$,}
  et  \vspace{0,1cm}\centerline{ $\ga(u_i,u_j)=\gve^{-1}.\ga(u_j,u_j)=\gve^{-1}.\go$.} 
 Mais la repr\'esentation $u$ \'etant non triviale, on a  $c\neq 0$, 
 donc  $\go\neq 0$, d'o\`u l'on d\'eduit que  $\gve=\gve^{-1}$ puis $|\gve|= |\gve_{i,j}|=1$ et enfin $|\go|= |c|$.
 \hfill $\Box$

{\bf \'Etude du cas particulier   $|\go|= |c|$}

On suppose encore que  la repr\'esentation $u : X\to (E,\ga)$ est r\'eduite, que  $|\go|= |c|$ et que $|X|\geq 3$.\\
Remarquons tout d'abord qu'on ne modifie pas le groupe des
automorphismes d'une gerbe $\GG(u)$ associ\'ee \`a une
repr\'esentation $u : X\to (E,\ga)$ en multipliant $u$ par un scalaire
non nul ou en rempla\cc ant $\ga$ par $-\ga$. L'\'etude du groupe
d'automorphismes de la gerbe $\GG(u)$ se ram\`ene donc toujours
 au cas o\`u $\go=1$, ce qu'on suppose dans la suite.
Quitte \`a r\'eordonner les \'el\'ements de $X$, on peut \'ecrire
$\GG(u)=\{\al u_1\ar\ldots,\al u_m \ar\}$ o\`u $m\leq n$ et les droites $\al
u_i\ar$ sont deux \`a deux distinctes pour $1\leq i\leq m$.
Notons $Y=\{1,\ldots,m\}$ et $(\gC_Y,Y)$ la restriction du graphe $(\gC,X)$ \`a $Y$, \cad le graphe dont l'ensemble des sommets est $Y$ et dont les ar\`etes sont celles de $(\gC,X)$ qui sont contenues dans $Y$. \\
Le \ttt suivant ram\`ene l'\'etude de la repr\'esentation $u : X\to
(E,\ga)$ et du groupe $\Aut(\GG(u))$ \`a celle de sa restriction $v$  \`a $Y$, qui rel\`eve du \ttt 2  :
\begin{thmf}$\ $\\ 
Soit  $v :Y\to (E,\ga)$ la restriction au graphe $(\gC_Y,Y)$ de la repr\'esentation $u :X\to (E,\ga)$. Alors,\\
$1.$ on a $\GG(u)=\{\al u_1\ar\ldots,\al u_n \ar\}=\{\al u_1\ar\ldots,\al u_m \ar\}=\GG(v)$\\
$2.$ les groupes $\Aut(\GG(u))$ et  $\Aut(\GG(v))$ sont \'egaux.\\
$3.$ la repr\'esentation $v :Y\to (E,\ga)$ du graphe $\gC_Y$ est
r\'eduite et non triviale.
\end{thmf}
\Dem. \\
 Les points 1 et 2 sont imm\'ediats.\\
 3.  L'\'egalit\'e  $\GG(u)=\GG(v)$ prouve que les matrices  $\Gram(u_1,\ldots,u_n)$ et $\Gram(v_1,\ldots,v_m)$ ont le \mm rang $r$, aussi \'egal \`a la dimension de $E$,  donc  $v :Y\to (E,\ga)$ est une repr\'esentation r\'eduite. Elle est non triviale car $u$ l'est.
\hfill $\Box$

{\it Commentaires}\\
En fait l'in\'egalit\'e stricte  $|\GG(u)|<|X|$ ne se produit que lorsque le graphe a une structure assez particuli\`ere que l'on regarde maintenant.\\
 {\it On introduit quelques notations.} 
Notons ''$\simeq$'' la relation \'equivalence sur $X$ donn\'ee par \ $i\simeq
k$ \ si \ $\al u_i\ar=\al u_k \ar$.  Notant  $Y=\{1,\ldots,m\}$,
chaque point $i$ de $X$ est donc \'equivalent  \`a un unique point $j$
de $Y$ que l'on note $j=\pi(i)$. Notons aussi $X_j=X_{\pi(i)}$ la
classe de l'\'el\'ement $i\in X$, et  pour $j\in Y$,   $X_j^+$
(resp. $X_j^-$)  d\'esigne l'ensemble des indices $i$ de $X_j$ tels
que $u_i=u_j$ (resp. $u_i=-u_j$). Enfin on applique aux exposants  la r\`egle des
signes usuelle, par exemple  $X_j^{+-}=X_j^-$.
\pb[1]
\begin{propf}$\ $\\ 
Soit $u$ une repr\'esentation du graphe $(\gC,X)$ de param\`etres $(1,c)$, o\`u $c=\pm 1$.\\
$1$. L'ensemble des ar\`etes qui  lient deux \elts  $A$ et $B$ parmi $X_1^+,X_1^-,\ldots, X_m^+,X_m^-$ est soit vide, soit l'ensemble de toutes les ar\`etes $\{a,b\}$ $($pour $(a,b)\in A\x B$$)$. \
On note  $A\not\sim B$ le premier cas $A\sim B$ le second.\\
$2.$ Pour deux indices $i$ et $j$ distincts dans $Y$ on a \\[0,1cm]
\vspace{0,1cm}\centerline{$X_i^+\sim X_j^+ \ssi X_i^+\not\sim X_j^-\ssi X_i^-\not\sim X_j^+\ssi X_i^-\sim X_j^-$}
$3.$ Pour tout indice $i$  dans $Y$ on a , \\[0,1cm]
\hspace*{0,4cm} $a$.\       si $c=1$,  \qquad $X_i^+\not\sim X_i^+$, \quad $X_i^-\not\sim X_i^-$\ et \  $X_i^+\sim X_i^-$, \\[0,1cm]
\hspace*{0,4cm} $b.$\    si $c=- 1$,  \quad $X_i^+\sim X_i^+$, \quad $X_i^-\sim X_i^-$\ et \  $X_i^+\not\sim X_i^-$.
\end{propf}
\Dem (la figure symbolique ci-dessous pr\'esente un schema possible du graphe dans les cas $c=1$ et $c=-1$ : chaque ar\`ete ou arc entre deux ensemble $A$ et $B$ symbolise le graphe form\'e de toutes les paires $\{a,b\}$ $($pour $(a,b)\in A\x B$$)$).\\
 \begin{center}
  \includegraphics[scale=0.25]{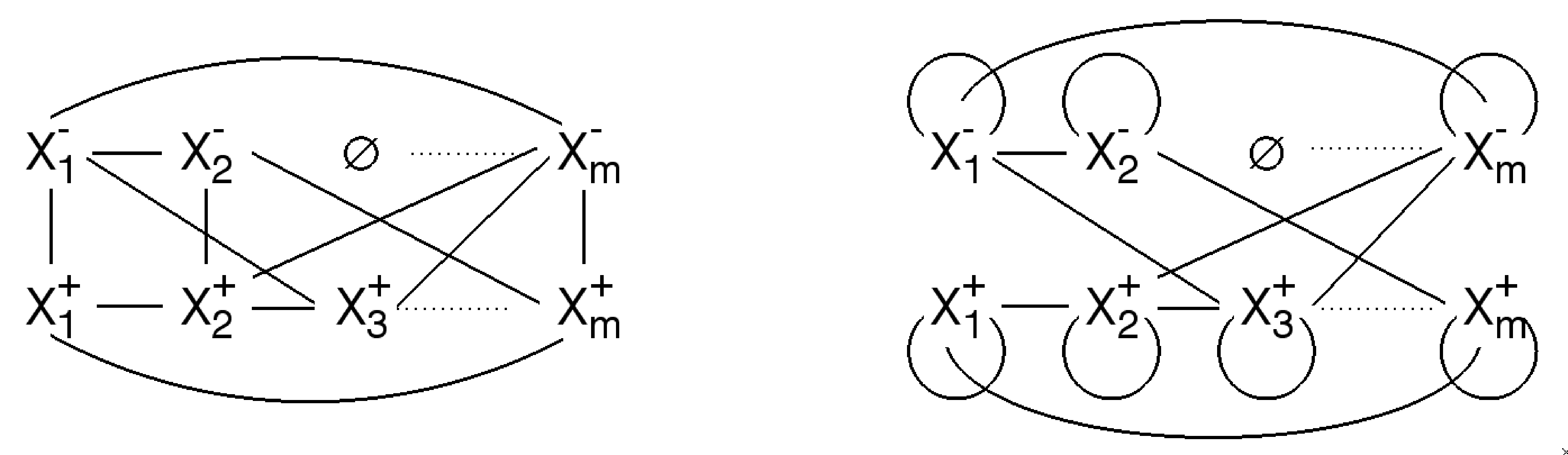}
   \npb \centerline{cas $c=1$\hspace{5cm} cas $c=-1$}
\end{center}
Soient  $A$ et $B$ deux ensembles parmi  $X_1^+,X_1^-,\ldots, X_m^+,X_m^-$, $i, i'$ deux indices pris dans $A$ et $j,j'$ deux indices pris dans $B$, donc  $u_{i'}=u_i$ et $u_{j'}= u_j$. \\
Si $A\neq B$ on a \\
\vspace{0,1cm}\centerline{
  $\gve_{i,j}.c=\ga(u_i,u_j)=\ga(u_{i'},u_{j'})=\gve_{i',j'}.c$,}
donc $\{i',j'\}$ est une ar\`ete \si $\{i,j\}$ en est une. De plus si $j''$ est dans $B^-$ alors \ $u_{j''}=-u_j$ \ et donc \\[0,1cm]
\vspace{0,1cm}\centerline{$\gve_{i,j''}.c=\ga(u_{j''},u_i)=-\ga(u_j,u_i)=-\gve_{i,j}.c$,}
 donc $\{i,j''\}$ est une ar\`ete \si $\{i,j\}$ n'en est pas une, ce qui  prouve 1 et 2 lorsque $A$ et $B$ sont distincts.\\
Si $A=B$ et $i\neq j$ comme  $u_i=u_j$, il vient $\gve_{i,j}.c=\ga(u_i,u_j)=\ga(u_i,u_i)=1$, donc le coefficient \ $\gve_{i,j}=c$ \ ne d\'epend que de $c$, et deux points de $A$ sont  li\'es si $c=-1$, et  non li\'es si $c=1$.
Ceci prouve 3 et termine la d\'emonstration. \hfill $\Box$

{\it Conclusion}\\
Les th\'eor\`emes  2 et 3 nous donnent une description compl\`ete du
groupe des automorphismes d'une gerbe $\GG(u)$ en fonction du graphe
$(\gC,X)$ dont elle provient.  Dans les exemples \'el\'ementaires qui
suivent, on montre que souvent le groupe $G(u)$ des automorphismes de
la gerbe est  plus gros que le groupe $H(\gC)$ des automorphismes du
graphe $(\gC,X)$.  Dans un travail en cours   nous utilisons ces
r\'esultats pour d\'eterminer toutes les gerbes isom\'etriques
$\GG(u)$ dont le groupe d'automorphismes agit deux fois transitivement sur $\GG(u)$.

\subsection{Quelques exemples}
Dans chacun des exemples qui suivent on donne :\\
Le graphe $(\gC, X)$ par sa matrice $\EE=S(1,1)$, la matrice  $S(1,c)$ et son d\'eterminant  $\chi(c ) = \det(S(1,c ))$ dont les racines fournissent les repr\'esentations de rang $r<n$ du graphe. Si $c$ est une racine de $\chi(c )$ de multiplicit\'e $\mu(c)$, la matrice $S(1,c)$ est de rang $r=n-\mu(c)$, donc associ\'ee \`a une repr\'esentation du graphe $(\gC,X)$ dans un espace de dimension $r$. On repr\'esente dans chaque cas un syst\`eme de vecteurs $u_1, \ldots, u_n$ qui r\'ealisent la repr\'esentation correspondante.
Pour des raisons purement graphique on a limit\'e les exemples aux cas o\`u la forme quadratique donn\'ee par la matrice $S(1,c)$ est positive.

\subsubsection{Triangle}
 $\EE=\left[ \begin {array}{ccc} 1&-1&-1\\\noalign{\medskip}-1&1&-1
\\\noalign{\medskip}-1&-1&1\end {array} \right] 
 $. \quad $S(1,c)=\left[ \begin {array}{ccc} 1&-c&-c\\\noalign{\medskip}-c&1&-c
\\\noalign{\medskip}-c&-c&1\end {array} \right] 
 $, \quad $\XX=- \left( 2\,c -1 \right)  \left( c +1 \right) ^{2} $ \\
  \begin{center}
  \includegraphics[scale=0.30]{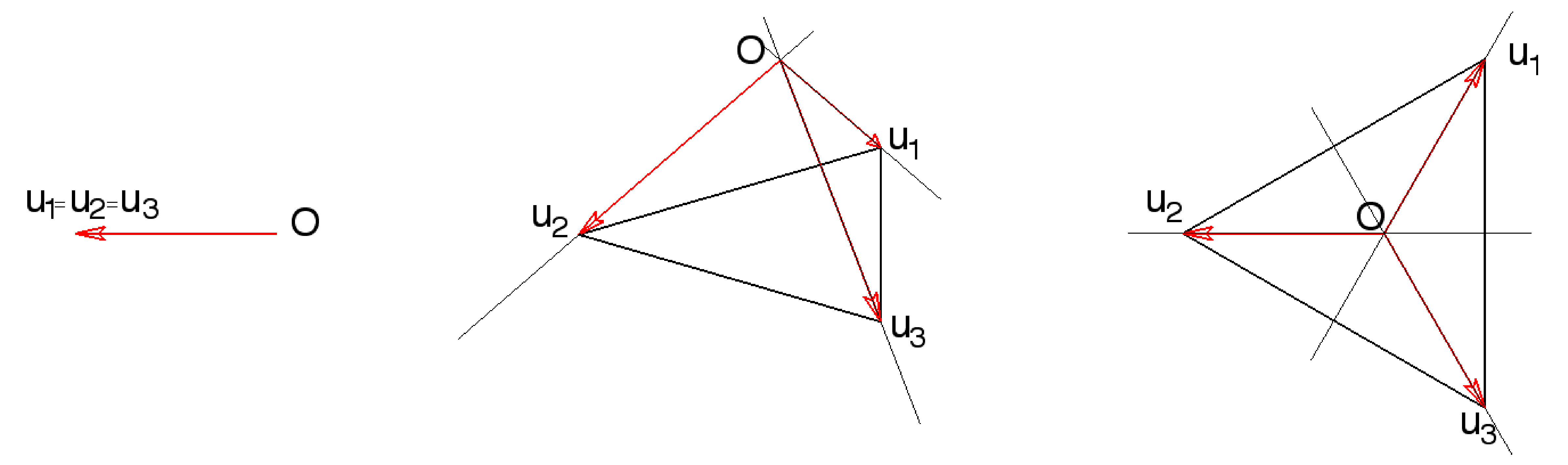}
   \npb \centerline{$c=-1$\hspace{3cm} $-1 <c  <1/2$\hspace{3,5cm}$c=1/2$}
\end{center}
Pour $c =-1$, la matrice $S(1,c )$ est de rang $1$ et  le triangle est repr\'esent\'e dans la droite $\R$. Les trois sommets du triangle sont envoy\'es sur un \mm point, par exemple $1$. Ce type de repr\'esentation (sur un point) ne peut avoir lieu que pour les graphes complets. \\
Pour $c =1/2$, la matrice $S(1,c )$ positive, est de rang $2$ et  le triangle est repr\'esent\'e dans le plan euclidien $\R^2$. C'est sa repr\'esentation usuelle comme triangle \'equilat\'eral, centr\'e \`a l'origine.\\
Pour $-1 <c  <1/2$, la matrice $S(1,c )$  est d\'efinie positive de rang $3$. Le triangle est repr\'esent\'e par les trois vecteurs issus d'un \mm sommet dans un t\'etra\`edre droit de base un triangle \'equilat\'eral. Une telle repr\'esentation est une combinaison lin\'eaire \`a coefficients positifs des deux pr\'ec\'edentes.\\
Lorsque $c  \not\in [-1,1/2]$ la repr\'esentation est de degr\'e trois, dans un espace quadratique donn\'e par une forme non positive.

\subsubsection{Carr\'e}
  $\EE= \left[ \begin {array}{cccc} 1&-1&1&-1\\\noalign{\medskip}-1&1&-1&1
\\\noalign{\medskip}1&-1&1&-1\\\noalign{\medskip}-1&1&-1&1\end {array}
 \right] 
$, \   $S(1,c)= \left[ \begin {array}{cccc} 1&-c&c&-c\\\noalign{\medskip}-c&1&-c&c
\\\noalign{\medskip}c&-c&1&-c\\\noalign{\medskip}-c&c&-c&1\end {array}
 \right] 
$, \   
$\XX=- \left( 3\,c +1 \right)  \left( c -1 \right) ^{3}$ \\
\begin{center}
  \includegraphics[scale=0.25]{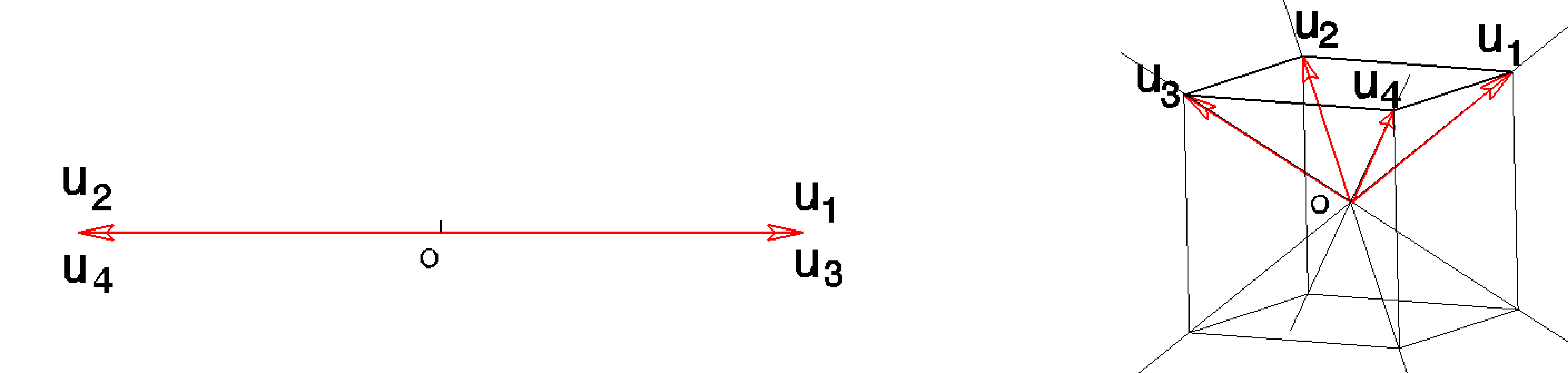}
   \npb \centerline{\hspace{2cm}$c=1$\hspace{7cm}$c=-1/3$}
\end{center}
Pour $c =-1/3$,  la matrice $S(1,c )$  est positive de rang $3$. Le carr\'e est repr\'esent\'e par la gerbe des quatre grandes diagonales d'un cube, qui font deux \`a deux un angle dont le cosinus vaut $-1/3$ comme pr\'evu. Le graphe du carr\'e est alors "concr\'etis\'e" par l'une des faces du cube. \\
Sur cet exemple le \ttt 2 prend tout son sens : le groupe d'automorphisme du carr\'e est le groupe di\'edral $D_4$ alors que le groupe de la gerbe associ\'ee est le groupe du cube, d'ordre $48$.\\
Pour $c = 1$,  la matrice $S(1,c )$  est positive de rang $1$. Le carr\'e est repr\'esent\'e sur la droite $\R$. Il s'envoie donc sur $\{-1,+1\}$. Deux sommets du carr\'e ont \mm image lorsqu'ils sont diagonalement oppos\'es.\\
Pour $-1/3 <c  <1$, la matrice $S(1,c )$  est d\'efinie positive de rang $4$. On visualise plus difficilement les gerbes associ\'ees toutefois ce sont toujours des combinaisons lin\'eaires \`a coefficients positifs des deux pr\'ec\'edentes. \\
Lorsque $c  \not\in [-1/3,1]$ la repr\'esentation est de degr\'e trois, dans un espace quadratique donn\'e par une forme non positive.

\subsubsection{Pentagone}
$\EE= \left[ \begin {array}{ccccc} 1&-1&1&1&-1\\\noalign{\medskip}-1&1&-1&1
&1\\\noalign{\medskip}1&-1&1&-1&1\\\noalign{\medskip}1&1&-1&1&-1
\\\noalign{\medskip}-1&1&1&-1&1\end {array} \right] 
$,  \  $S(1,c)= \left[ \begin {array}{ccccc} 1&-c&c&c&-c\\\noalign{\medskip}-c&1&-c&c
&c\\\noalign{\medskip}c&-c&1&-c&c\\\noalign{\medskip}c&c&-c&1&-c
\\\noalign{\medskip}-c&c&c&-c&1\end {array} \right] $, 
 \  $\XX=\left( -1+5\,{c }^{2} \right) ^{2} $ \\[0,3cm]
Pour $c = \pm\sqrt{5}/5$,   la matrice $S(1,c )$  est positive de
rang $3$.  On obtient ces repr\'esentations en consid\'erant cinq des  six droites passant par les centres des faces oppos\'ees d'un dod\'eca\`edre. Elles  forment deux \`a deux un angle dont le cosinus est $\pm \sqrt{5}/5$. 
Si $c=\sqrt{5}/5$,  pour $j=i+1 \mod(5)$, les points $i$ et $j$ de $X=\{1,2,3,4,5\}$ sont li\'es, donc  $(u_i|u_j)=-\sqrt{5}/5$. On obtient la repr\'esentation crois\'ee du pentagone, tandis que pour $c=-\sqrt{5}/5$,  les m\^emes points $i$ et $j$ sont li\'es, donc  $(u_i|u_j)=\sqrt{5}/5$ ce qui donne la repr\'esentation '' convexe'' du pentagone.\\
Le groupe d'automorphismes de la gerbe des 5 droites $U_1=\al u_1\ar,\ldots,U_5=\al u_5\ar$ est le stabilisateur de $U_6$ dans le groupe des isom\'etries de l'icosa\`edre. Il est donc isomorphe au groupe $D_5$.\\
 \begin{center}
  \includegraphics[scale=0.25]{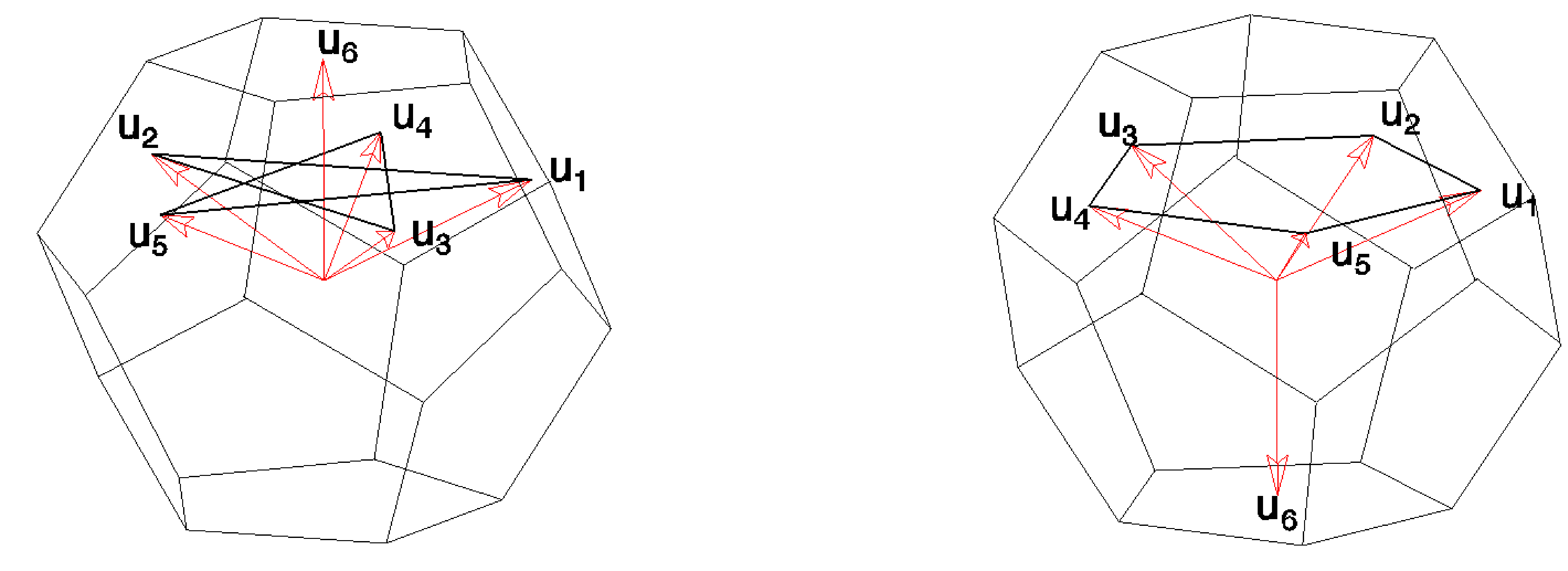}
\npb \centerline{$c = \sqrt{5}/5$\hspace{6cm} $c = -\sqrt{5}/5$}
 \end{center}

\subsubsection{Hexagone point\'e}
J'appelle hexagone point\'e le graphe  form\'e sur six points par un pentagone et un point isol\'e. Ce cas est int\'eressant parce qu'il sera fortement g\'en\'eralis\'e par la suite.\\[0,2cm]
$\EE=\left[ \begin {array}{cccccc} 1&-1&1&1&-1&1\\\noalign{\medskip}-1&1&-
1&1&1&1\\\noalign{\medskip}1&-1&1&-1&1&1\\\noalign{\medskip}1&1&-1&1&-
1&1\\\noalign{\medskip}-1&1&1&-1&1&1\\\noalign{\medskip}1&1&1&1&1&1
\end {array} \right] 
 $, \ $S(1,c)=\left[ \begin {array}{cccccc} 1&-c&c&c&-c&c\\\noalign{\medskip}-c&1&-
c&c&c&c\\\noalign{\medskip}c&-c&1&-c&c&c\\\noalign{\medskip}c&c&-c&1&-
c&c\\\noalign{\medskip}-c&c&c&-c&1&c\\\noalign{\medskip}c&c&c&c&c&1
\end {array} \right] 
 $,
  \ 
$\XX=- \left( -1+5\,{c }^{2} \right) ^{3}$ \\[0,2cm]
Comme on le voit le polyn\^ome $\chi(c )$ est un multiple du pr\'ec\'edent, ce qui n'est gu\`ere \'etonnant puisque le pentagone apparait comme un sous-graphe de l'hexagone point\'e. \\
Les deux repr\'esentations, associ\'ees aux cas $c = \sqrt{5}/5$ et $c = -\sqrt{5}/5$ proviennent des repr\'esentations  du pentagone auxquelles on rajoute le vecteur $u_6$. Mais comme le point $6$ n'est li\'e \`a aucun des points $1,2,3,4,5$, l'angle $(u_i,u_6)$ (pour $1\leq i \leq 5$) est le compl\'ementaire des angles $(u_i,u_{i+1})$ ($i$, entier modulo  $5$), autrement dit si l'angle $(u_i,u_{i+1})$ est aigu alors $(u_i,u_6)$ est obtus et vice versa. 
Dans les deux cas le groupe $H$ d'automorphismes du graphe est le groupe di\'edral $D_5$ tandis que le groupe d'automorphismes de la gerbe image est le groupe d'isom\'etries du dod\'eca\`edre, agissant deux fois transitivement sur l'ensemble $\GG(u)=\{U_1,\ldots , U_6\}$ des six droites vectorielles.

\section{Bibliographie}

\label{thebib}
{\it Livres}\\[0,2cm]
$[1]$\quad  Dembowski, P. {\it   Finite  geometries.} Springer-Verlag (1968)\\[0,1cm]
$[2]$\quad  Norman Biggs. {\it Finite groups of automorphisms.} \mbox{London Mathematical Society. (1970)}\\[0,1cm]
$[3]$\quad  Norman Biggs. {\it Algebraic Graph Theory.} \mbox{Cambridge University Press,  (1993)}\\[0,1cm]
$[4]$\quad P. J. Cameron and J.H.Van Lint. {\it Designs, graphs,codes and their links, vol 22 of London Mathematical Society Student Texts,} Cambridge University Press, Cambridge, 1991.
\end{document}